\documentclass[14pt]{amsart}
\usepackage{amssymb}
\usepackage{enumerate}
\usepackage[colorlinks]{hyperref}

\newtheorem{thm}{Theorem}[section]

\newtheorem{thm-def}[thm]{Theorem-Definition}
\newtheorem{lem}[thm]{Lemma}

\theoremstyle{definition}

\theoremstyle{remark}
\newtheorem{rem}[thm]{Remark}
\numberwithin{equation}{section}

\DeclareMathOperator{\tensor}{\otimes}

\newcommand{\vacB}{|0\rangle}
\newcommand{\vac}{\vacB}

\newcommand{\nc}{\newcommand}

\nc{\DMO}{\DeclareMathOperator}

\DMO{\Spec}{Spec}

\nc{\Real}{\mathbb R}
\nc{\Cplx}{\mathbb C}
\nc{\ZZ}{\mathbb Z}
\nc{\LL}{\mathbb L}
\nc{\VV}{\mathbb V}
\nc{\MM}{\mathbb M}

\nc{\Field}{\mathbb F}
\nc{\Rat}{\mathbb Q}

\nc{\gr}{{\rm gr\ }}
\nc{\grG}{{\rm gr^{\mathcal{G}}\ }}
\nc{\grF}{{\rm gr^{\mathcal{F}}\ }}
\nc{\grH}{{\rm gr^{\mathcal{H}}\ }}

\nc{\pole}{\Cplx}

\nc{\de}{\partial}
\nc{\oX}{{\mathcal{O}}_{X}}

\nc{\calO}{{\mathcal{O}}}

\nc{\calA}{{\mathcal{A}}}
\nc{\calB}{{\mathcal{B}}}
\nc{\calC}{{\mathcal{C}}}
\nc{\calD}{{\mathcal{D}}}
\nc{\calE}{{\mathcal{E}}}

\nc{\calF}{{\mathcal{F}}}
\nc{\calG}{{\mathcal{G}}}
\nc{\calH}{{\mathcal{H}}}

\nc{\calK}{{\mathcal{K}}}

\nc{\calT}{{\mathcal{T}}}

\nc{\cA}{\mathcal{A}}
\nc{\cB}{\mathcal{B}}

\nc{\cD}{{\mathcal{D}}}

\nc{\cF}{{\mathcal{F}}}
\nc{\cG}{{\mathcal{G}}}
\nc{\cH}{{\mathcal{H}}}
\nc{\cJ}{{\mathcal{J}}}

\nc{\cL}{{\mathcal{L}}}
\nc{\cM}{{\mathcal{M}}}
\nc{\cN}{\mathcal{N}}
\nc{\cO}{{\mathcal{O}}}
\nc{\cQ}{{\mathcal{Q}}}

\nc{\cS}{{\mathcal{S}}}

\nc{\cT}{{\mathcal{T}}}

\nc{\cV}{{\mathcal{V}}}

\nc{\cZ}{{\mathcal{Z}}}

\nc{\TT}{\partial} 

\nc{\wt}[1]{\Delta_{#1}}

\nc{\op}[1]{{#1}_{(-1)} }

\nc{\opp}[2]{ {#1}_{(  #2  )} }

\nc{\zero}{{}_{(0)} }
\nc{\one}{{}_{(1)} }
\nc{\opm}{{}_{(-1)} }
\nc{\ops}[1]{{}_{({#1})} }

\DMO{\Hom}{Hom}
\DMO{\Aut}{Aut}
\nc{\Endom}{{\rm End\,}}
\nc{\End}{{\mathrm End\,}}

\nc{\Iso}{{\rm Iso}}

\DeclareMathOperator{\id}{id}

\nc{\Ind}{{\rm Ind\ }}
\nc{\im}{{\rm Im\ }}

\DMO{\spn}{span}
\DMO{\Der}{Der}

  \nc{\be}{\begin{equation}}
  \nc{\ee}{\end{equation}}

\nc{\iso}{\simeq}

\nc{\res}{{\rm res}}
\nc{\Resz}{{\rm Res}_z}
\nc{\ta}{{\tilde{a}}}
\nc{\tb}{{\tilde{b}}}
\nc{\ttt}{{\tilde{t}}}
\nc{\ts}{{\tilde{s}}}
\nc{\tg}{{\tilde{g}}}
\nc{\tf}{{\tilde{f}}}

\nc{\txi}{{\tilde{\xi}}}
 \nc{\teta}{{\tilde{\eta}}}

\nc{\tlam}{{\tilde{\lambda}}}
	
\nc{\la}{\lambda}
\nc{\La}{\Lambda}

\nc{\fg}{\frak{g}}
\nc{\fb}{\frak{b}}
\nc{\fn}{\frak{n}}
\nc{\fh}{\frak{h}}
\nc{\fz}{\frak{z}}
\nc{\fzl}{\fz^{\text{reg},\lambda}}

\nc{\fU}{{\frak{U}}}
\nc{\fM}{{\frak{M}}}

\nc{\slt}{\frak{s}\frak{l}_2}
\nc{\affsl}{\widehat{\frak{s}\frak{l}}_2}

\nc{\prline}{\mathbb{P}^1}

\nc{\Vcrit}{V_{\kappa_c}(\slt)}
\nc{\Symb}{{\rm Symb\  }}

\nc{\omcl}{\Omega_X^{1,\,cl}}
\nc{\homcl}{H^1 (X, \omcl)}
\nc{\funomega}{Fun(\homcl)}

\nc{\cdo}{{\cD}_X^{ch}}
\nc{\dl}{\cD^\cL}
\nc{\twcdo}{{\cD}_X^{ch, tw}}
\nc{\twcdoloc}{{\stackrel{\circ}{\cD}}_X^{ch, tw}}
\nc{\twcdolocflag}{{\stackrel{\circ}{\cD}}_{G/B}^{ch, tw}}
\nc{\gtdo}{{\cD}_X^{tw}}
\nc{\cdofun}[1]{\cD_{#1}\otimes\funomega}

\nc{\Pic}{ {\rm Pic}  }
\nc{\cch}{\Cplx_{\chi(z)}}
\nc{\cps}{\Cplx_{\Psi}}
\nc{\dchmod}{\cdo{\textrm{-}}Mod }
\nc{\dchtwmod}{\twcdo{\textrm{-}}Mod }
\nc{\dmod}{\cD_X{\textrm{-}}Mod}
\nc{\dlmod}{\cD^{\cL}_X{\textrm{-}}Mod}
\nc{\dlm}{\cD^{\cL}_X{\textrm{-}}Mod}

\nc{\dlamod}{\gtdo{\textrm{-}}Mod}
\nc{\tensL}{\otimes_{\oX}\!\cL}
\nc{\ttens}{\overline{\otimes}}
\nc{\bM}{\bar{M}}
\nc{\bg}{\bar{g}}
\nc{\btau}{\bar{\tau}}
\nc{\bxi}{\bar{\xi}}
\nc{\ba}{\bar{a}}

\nc{\abs}[1]{\left\vert#1\right\vert}
\nc{\set}[1]{\left\{#1\right\}}
\nc{\eps}{\varepsilon}
\nc{\To}{\longrightarrow}

\nc{\CDOP}{\cD^{ch}_{\prline}}
\nc{\TWCDOP}{\cD^{ch,tw}_{\prline}}
\nc{\twcdop}{\cD^{ch,tw}_{\prline}}

\nc{\basis}{ \set{\tau^i} }
\nc{\basisom}{\{ \omega_i\}}

\nc{\Lie}{{\rm Lie\,}}
\nc{\lie}[1]{ {\, \rm Lie}^{}_{\,#1}\,  }

\nc{\cfdeg}{{\rm conf.deg}}

\nc{\PP}{P(\de, \omega)}
\nc{\tP}{\tilde{P}}
\nc{\tQ}{\tilde{Q}}

\nc{\partialtau}[1]{D_{\de_{#1}}}

\nc{\partialom}[1]{D_{\omega_{#1}}}

\nc{\Fields}{{Fields}}
\nc{\ghat}{{\hat{\fg}}}
\nc{\qf}[2]{\langle {#1}, {#2}  \rangle}
\nc{\vlk}{v_{\lambda, k}}
\nc{\uH}{\underline{\rm{H}}}

\nc{\pair}[1]{\langle #1 \rangle}

\nc{\homeg}{ H^1 (X, \Omega_X^{ [1, 2 \rangle} ) }

\nc{\ttau}{ {\tilde{\tau}} }
\nc{\tla}{{\tilde{\lambda}^*}}

\nc{\tfg}{ {\tilde{\fg}} }

\nc{\itaul}{{i^{}_{\tau_l} }}

\nc{\itaum}{{i^{}_{\tau_m} }}

\nc{\ixi}{{i^{}_\xi}}
\nc{\ieta}{{i^{}_\eta}}
\nc{\half}{\frac{1}{2} }

\nc{\tPhi}{\tilde{\Phi}}
\nc{\gV}{{\fg(\cV) }}

\nc{\disp} {\displaystyle}

\DeclareMathOperator{\Sing}{Sing}

\nc{\Gr}{{\mathcal{G}r}}

\nc{\CDO}{{\mathcal{C}\mathcal{D} \cO}}

\nc{\Om}[3]{{ \Omega_{#3}^{#1} \to \Omega_{#3}^{#2, cl}}}

\nc{\byzero}{ \left[ \Sing M \right]}

\nc{\oinv}{{ M^{\Omega} }}
\nc{\dd}{{\mathbf{d}}}
 
\nc{\Moh}{M^{\Omega, \fh}}

\nc{\cz}{{\mathfrak{z}}}

\nc{\fhn}{\fh^{\nabla}}
\nc{\halfint}{\cM^{int_+}}
\nc{\halfinttheta}{\halfint_{\theta}}
\nc{\halfintchitheta}{\halfint_{\theta, \chi(z)} }
\nc{\thetabar}{\bar{\theta}}

\begin{document}
\pagestyle{plain}
\title{Half-integrable modules  over algebras of twisted chiral differential operators}
\author{D. Chebotarov}
\maketitle
\begin{abstract}
A module $M$ over a vertex algebra $V$ is {\em half-integrable}
 if $a_n$ act locally nilpotently on $M$ for all  $a\in V$, $m\in M$, $n>0$
 We study half-integrable modules over sheaves of twisted chiral differential operators (TCDO) 
 on a smooth variety $X$. 
 We prove an equivalence  between certain categories of half-integrable modules over TCDO
  and categories of (twisted) D-modules on $X$.
\end{abstract}

 \section{Introduction}
 
 Let $V$ be a graded vertex algebra and $M$ a $V$-module, i.e. we are given a map
 $$
 V_\Delta \ni a \mapsto a^M (z) = \sum a_n z^{-n-\Delta} \in End M [[z,z^{-1} ]]
 $$
 satisfying the usual axioms, cf. section \ref{Definitions} below.
  
 We call  $M$ {\em half-integrable}
 if for any $a \in V$ and $n>0$ 
 the action of
 $a_n$ on $M$ is locally nilpotent; 
that is, $M$ is integrable with respect to the positive part of the graded Borcherds Lie algebra of $V$.

 Let $X$ be a  smooth complex variety
 and $\cV$  a graded sheaf of vertex algebras on $X$
 which is a vertex envelope of   vertex algebroid $\cV_0 \oplus \cV_1$, see \cite{GMS}, such that  $\cV_0 = \cO_X$. 
 The previous definition immediately carries over to this situation, and we will apply it to those $\cV$ that satisfy some extra conditions.
 In order to formulate these conditions, recall that 
$\cL := \cV_1 / \cV_0 \opm \de \cV_0$ is an $\cO_X$-Lie algebroid \cite{GMS};
in particular, there is an anchor, an $O_X$-Lie algebroid
morphism: $\cL \to \cT_X$. 
We require:
\begin{itemize}
 \item 
 $\cL$ is a  locally free $\cO_X$-module of finite rank;
  \item
     $\cL$ fits into an exact sequence 
     $$
     0 \To \fh   \To \cL  \To   \cT_X   \To 0
     $$
     where $\fh$ is an abelian $\cO_X$-Lie algebra.
  \end{itemize}
  Such sheaves of vertex algebras $\cV$ form a natural class which includes  sheaves of chiral twisted differential operators $\cD^{ch, tw}_X$ defined in \cite{AChM} 
  and certain deformations of those. 
	  
  \smallskip

 Let $\cZ (\cV) \subset \cV$ be the subsheaf of $\cV$ 
 with stalks $\cZ_{X,x} = Z(\cV_{X,x})$, $x\in X$,
 where $Z(V)$ denotes the center of $V$; we will call it 
 the {\em center} of $\cV$.

 Define $\halfint(\cV)$ to be the category of half-integrable  $\cV$-modules 
 that satisfy the following regularity condition:
  $z_n m=0$ for $z\in \cZ(\cV)$, $n>0$, $m\in \cM$.

 Then the following result  holds. (Theorem \ref{main_result_gen}).
 
 \smallskip

{\bf Theorem 1.} 
{ \em 
Let $\cV$ be as above and $\cM \in \halfint(\cV)$. Then

1) \, The subsheaf 
 $$
  \Sing \cM = \set{ m \in \cM : \ v_n m = 0  \textrm{ for all } v\in \cV,  \,  n>0 }
 $$
   generates $\cM$;
   
   2) $\cM$ possesses a filtration 
    $\cM = \cup_{i\ge 0} \cM_i$ 
    such that 
  $\cV_{i} {}_{(n)}  \cM_{j}  \subset \cM_{i+j - n -1}$
  and $\cM_0 = \Sing \cM$.
   } 
 
 \smallskip

The result, as stated, is true in the analytic topology; 
when working algebraically one has to add one more requirement,
cf. section \ref{the-algebraic-case-Section}.

\bigskip

  Let $X$ be connected and $\cV$ as above.
 Suppose that
 $\fh = \cO_X \tensor_\Cplx \fhn$
where 
  $\fhn$ is a constant sheaf such that $[\cL, \fhn] = 0$.
 For 
 $\theta \in (\fhn)^*$  
 consider the ideal $I_\theta$ in $U_{\cO_X} (\cL)$ 
 generated by $h - \theta(h)$, $h \in \fhn$.
 Then 
 $U_{\cO_X} (\cL) / I_\theta$ 
 is an algebra of twisted differential operators (tdo) on $X$. 
 Let us denote it by 
 $\cD_X^{\thetabar}$.

 Let $\halfinttheta(\cV)$ denote the full subcategory of modules
  $\cM\in \halfint(\cV)$ such that
 for each  $h \in \fhn$, $m \in \cM$ \  $h_0 m = \theta(h) m$.
 Being the Zhu algebra of $\cV$,
  $U_{\cO_X} (\cL)$ acts on $\cM_0 = \Sing \cM$,
  and this action factors through   
 $U_{\cO_X} (\cL) / I_\theta = \cD_X^{\thetabar}$. 
 This way we get a functor
 \begin{equation}
 \label{intro-sing-hinttheta}
  \Sing:  \halfinttheta(\cV)    \to   \cM(\cD_X^{\thetabar})
 \end{equation}
 It is not an equivalence;  however, it becomes one upon restriction to certain subcategories. 

 It is not hard to verify that the center $\cZ(\cV)$ is a constant sheaf 
 with values in an algebra 
  of differential polynomials, generated by the subspace
 $\cz: = (\fh^{\nabla})^{\perp} = \set{ h\in \fh^{\nabla} \, :\ \pair{h, \fh^{\nabla}} = 0 }$, 
 where $\pair{ , }: \fh \times \fh \to \cO_X$ is induced by $\one$.
 
  Fix a series $\chi(z) \in \cz^*((z))$,  
  $\chi(z)  = \sum_{}  \chi_n z^{-n-1}$. 
  We say that $\cZ$ acts on $\cM$ via $\chi(z)$ if for all $h\in \cz$, $m\in \cM$
  $$
  h_n m = \chi_n (h) m, \ n\in \ZZ
  $$ 
 Let $\halfintchitheta(\cV)$ denote the  full subcategory of $\halfinttheta(\cV)$ consisting of modules with action of $\cZ$ given by $\chi(z)$. 
 For this category to be nonzero, the obvious compatibility conditions 
 $\chi_0 = \theta|_{\cz}$ and $\chi(z) \in \cz^*[[z]]z^{-1}$
 have to be satisfied. 
 
 The following result holds. (Theorem \ref{equivalence-theorem} below).
 
 \smallskip
 
 {\bf Theorem 2.}  
 {    \em
  Let $\chi(z)\in \cz^*[[z]]$ be such that  $\chi_0 = \theta|_{\cz}$.
    Then the  restriction of (\ref{intro-sing-hinttheta})   
    to $\halfintchitheta(\cV)$ is an equivalence.
   }
   
   \smallskip

 As a particular case, this theorem can be applied when $\cV$ is the sheaf of twisted chiral differential operators $\cD^{ch, tw}_X$ introduced in \cite{AChM}. In this case $\fhn$ and $\cz$ are both equal to  the trivial local system with fiber $H^1 (X, \Omega^{1} \to \Omega^{2,cl})^*$ 
 and the center $\cZ(\cV)$ is the algebra of differential polynomials on the affine space
 $H^1 (X, \Omega^{1} \to \Omega^{2,cl})$. 
Theorem 2  becomes an equivalence between 
 the category of half-integrable $\twcdo$-modules with central character $\chi(z)$  and 
 the category of modules over the tdo $\cD^{\chi_0}$. 
 
 Furthermore, when $\fh = 0$
 this becomes an equivalence between the category of half-integrable 
  modules over a CDO $\cD^{ch}_X$ and the category of usual $\cD_X$-modules.

These results should be contrasted with the classic Kashiwara Lemma from the
D-module theory (\cite{Borel}, \cite{HTT}). Recall
that if $Y\subset X$ is a submanifold and $J$ is the defining ideal,
$M$ is a $D_X$-module supported on Y, and $M_J\subset M$ the subsheaf annihilated
by $J$, then $M_J$ is naturally a $D_Y$-module. Kashiwara's Lemma asserts
that $M \mapsto M_J$ sets up an equivalence of the category of $D_X$-modules
supported on $Y$ and the category of $D_Y$-modules.

A    TCDO can be intuited about as an algebra of differential operators
on the loop space $\cL M$, see papers of Kapranov and Vasserot (\cite{KV1, KV2}) for
a rigorous treatment. The half-integrability condition can then be interpreted as the requirement that the module is supported on regular loops.
$\Sing M$ then becomes an analogue of $M_J$.
Therefore, Theorem 2 appears to be a loop space (or vertex algebra) version  of Kashiwara's Lemma.

{\em Acknowledgement.} 
The author is  deeply grateful to Fyodor Malikov for invaluable help 
throughout the work.

\bigskip

\section{Preliminaries.}

We will recall the basic  notions of vertex algebra 
following the exposition of \cite{AChM}.

All vector spaces will be over  $\pole$. 
All vertex algebras considered in this article will be even  $\ZZ_{\ge 0}$-graded vertex algebras.

\subsection{Definitions.}
\label{Definitions}
 Let $V$ be a vector space.

A {\em field} on $V$ is a formal series $$a(z) = \sum_{n\in \ZZ}
a_{(n)} z^{-n-1} \in ({\rm End} V)[[z, z^{-1}]]$$ such that for any
$v\in V$ one has $a_{(n)}v = 0$ for sufficiently large $n$.

Let $\Fields (V)$ denote the space of all fields on $V$.

A {\em vertex algebra}  is a vector space $V$ with the following
data:
\begin{itemize}
  \item a linear map $Y: V \to \Fields(V)$,  
    $V\ni a \mapsto a(z) = \sum_{n\in \ZZ} a_{(n)} z^{-n-1}$
  \item a vector $\vac\in V$, called {\em vacuum vector}
  \item a linear operator $\de: V \to V$, called {\em translation operator}
\end{itemize}
that satisfy the following axioms:
\begin{enumerate}
  \item (Translation Covariance)

 $ (\de a)(z) = \de_z a(z)$

  \item (Vacuum)

  $\vac(z) = \id$;

  $a(z)\vac \in V[z]$ and $a_{(-1)}\vac = a$

  \item (Borcherds identity)
   \begin{align}
   \label{Borcherds-identity}
   &\sum\limits_{j\geq 0} {m \choose j} (a\ops{n+j} b )\ops{m+k-j}\\
   =& \sum\limits_{j\geq 0} (-1)^{j} {n \choose j}\{ a\ops{m+n-j} b\ops{k+j} - (-1)^{n}b\ops{n+k-j} a\ops{m+j}
   \}\nonumber
   \end{align}
\end{enumerate}

A vertex algebra $V$ is {\em graded} if  $V = \oplus_{n\geq 0}V_n$
and for $a\in V_i$, $b\in V_j$ we have
$$a_{(k)}b \in V_{i+j - k -1}$$ for all $k\in \ZZ$. (We put $V_i  = 0$ for $i<0$.)

  All vertex algebras in this article are graded vertex algebras.

We say that a vector $v\in V_m$ has {\em conformal weight} $m$ and
write $\Delta_v = m$.

If $v\in V_m $ we denote $v_k  = v_{(k - m +1)}$, this is the
so-called conformal weight notation for operators. One has 
$$v_k V_m
\subset V_{m -k}.$$

 A {\em morphism} of vertex algebras is a map $f: V \to W$ that preserves vacuum and satisfies $f(v_{(n)}v') = f(v)_{(n)}f(v')$.

A {\em module} over a vertex algebra $V$ is a vector space $M$
together with a map
\begin{equation}
\label{def-vert-mod-1} Y^M: V \to \Fields(M),\;  a \to Y^M(a,z) =
\sum_{n\in \ZZ} a^M_{(n)}z^{-n-1},
\end{equation}
 that satisfy the following axioms:
\begin{enumerate}
  \item $\vac^M (z)   = \id_M  $
  \item (Borcherds identity)
  \begin{align}\label{def-vert-mod-2}
     &\sum\limits_{j\geq 0} {m \choose j} (a^{}_{(n+j)} b
     )^M_{(m+k-j)}\\
  = &\sum\limits_{j\geq 0} (-1)^{j} {n \choose j}\{ a^M_{(m+n-j)} b^M_{(k+j)} - (-1)^{n}b^M_{(n+k-j)} a^M_{(m+j)} \}\nonumber
  \end{align}
\end{enumerate}

\smallskip

A module $M$ over a graded vertex algebra $V$ is called {\em graded}
if $M = \oplus_{n\geq 0} M_n$ with
 $v_{k}M_l  \subset M_{l-k}$  (assuming $M_{n} = 0$ for negative $n$).

 A {\em morphism of modules} over a vertex algebra $V$ is a map $f: M \to N$
 that satisfies $f(v^M_{(n)}m) = v^N_{(n)}f(m)$ for $v\in V$, $m\in M$.
$f$ is {\em homogeneous} if $f(M_k)\subset N_k$ for all $k$.

\subsubsection{The Borcherds Lie algebra}
\label{BorcherdsLie}
To any vertex algebra $V$ one can associate 
a Lie algebra that acts on any $V$-module. It is 
 the {\em Borcherds Lie algebra} of $V$ defined in the following way
$$
\Lie (V) = V\tensor \Cplx[t,t^{-1}]  \,  /  \, (\,(\de a  + (n + H) a)\tensor t^n) 
$$
where $Ha = ka$ for $a\in V_k$.
This is a Lie algebra with Lie bracket given by
$$
[a\tensor t^n, b\tensor t^l] = \sum_{j\geq 0} { n +\Delta_a -1 \choose j} (a\ops{j} b)\tensor t^{n+l} 
$$
 for a homogeneous $a\in V_\Delta$ and $b\in V$,
and extended linearly.
   $\Lie(V)$  acts on any $V$-module $M$ by letting $a\tensor t^n$ act as $a_n$.

	 $\Lie(V)$  has a natural grading, $\disp \Lie(V) = \oplus_{n \in \ZZ} \, \Lie(V)_n $  with $\Lie(V)_n$ equal to the image of $V\tensor t^n$.

\smallskip

\subsection{Examples}

\subsubsection{Commutative vertex algebras.}
\label{Commutative_vertex_algebras} A vertex algebra is said to be
  {\em commutative} if $a_{(n)} b =0$ for $a$, $b$ in $V$ and $n\geq 0$.
   The structure of a commutative vertex algebras is equivalent to one 
 of commutative associative algebra with a derivation.

If $W$ is a vector space we denote by $H_W$ the algebra of
differential polynomials on $W$. As an associative algebra it is a
polynomial algebra in variables $x_i$, $\de x_i$, $\de^{(2)}x_i$,
$\dots$ where $\set{x_i}$ is a basis of $W^*$. A commutative vertex
algebra structure on $H_W$ is uniquely determined by attaching the
field
 $x(z) = e^{z\de}x_i$  to $x\in W^*$.

 $H_W$ is equipped with grading such that
 \begin{equation}
 \label{weights-0-1-diff-poly}
 (H_W)_0=\pole,\; (H_W)_1=W^*.
 \end{equation}

\subsubsection{Beta--gamma system.}
\label{beta-gamma-system} Define the Heisenberg Lie algebra  to be
the algebra with generators $a^i_n$, $b^i_n$, $1\leq i\leq N$ and
$K$ that satisfy $[a^i_m, b^j_n] = \delta_{m, -n} \delta_{i,j} K$,
$[a^i_n, a^j_m] = 0$, $[b^i_n, b^j_m] = 0$.

Its Fock representation $M$ is defined to be the module induced from
the  one-dimensional representation $\Cplx_1$ of its subalgebra
spanned by $a^i_n$, $n\geq 0$, $b^i_m$, $m>0$ and $K$ with $K$
acting as identity and all the other generators acting as zero.

The {\em beta-gamma system} has $M$ as an underlying vector space,
the
 vertex algebra structure being determined  by
assigning the fields
$$a^i(z) = \sum a^i_n z^{-n-1}, \ \ b^i(z) = \sum b^i_n z^{-n}$$
to $a^i_{-1}1$ and $b^i_01$ resp., where $1\in \Cplx_1$.

This vertex algebra is given a grading so that the degree of
operators $a^i_n$ and $b^i_n$ is $n$. In particular,
\begin{equation}
\label{weight-0-1-b-g} M_0=\pole[b_0^1,...,b_0^N],\;
M_1=\bigoplus_{j=1}^{N}(b^j_{-1}M_0\oplus  a^j_{-1}M_0).
\end{equation}

\subsection{Vertex algebroids}
\label{Vertexalgebroids}

\subsubsection{}
\label{def-of-vert-alg}
	Let  $V$ be a  vertex algebra. 

	Define a 1-truncated vertex algebra to be a sextuple
	 $(V_0\oplus V_1, \vac, \de, \opm, \zero, \one )$ 
	 where the operations $\opm, \zero, \one$ satisfy all the axioms
	   of a vertex algebra that make sense upon restricting to the subspace $V_0 + V_1$.
	    (The precise definition can be found in \cite{GMS}).
	   The category  of 1-truncated vertex algebras will be denoted $\cV ert_{\le 1}$.

The definition of vertex algebroid is a 
reformulation of that of a sheaf of 1-truncated vertex algebras.

A {\em vertex  $\cO_X$-algebroid} is a sheaf $\cA$ of  $\Cplx$-vector spaces equipped with $\Cplx$-linear maps 
$\pi: \cA \to \cT_X$ and $\de: \cO_X \to \cA$ satisfying $\pi \circ \de = 0$
and with    
operations
$
\opm :  \cO_X \times \cA   \To \cA
$,
$
\zero : \,  \cA \times \cA    \To \cA
$,
$
\one :  \,  \cA \times \cA  \To \cO_X
$
satisfying axioms:
\begin{eqnarray}
f\opm(g\opm v) - (fg)\opm v & = & \pi(v)(f)\opm \partial(g) +
\pi(v)(g)\opm \partial(f)
\label{v-assoc} 
\\
 x \zero ( f\opm y )  & = & \pi(x)(f)\opm y + f\opm ( x \zero y )
 \label{leib}
\\
 x \zero y  +  y \zero x & = & \partial (  x \one y )
\label{symm-zero}
\\
\pi(f\opm v) & = & f\pi(v) 
\\
 ( f\opm x) \one y   & = & f ( x\opm y)  -
\pi(x)(\pi(y)(f)) 
\label{opm-one}
\\
 \label{zero-one} 
\pi(v)( x \one y ) & = & (  v \zero x) \one y  +
x \one (v \zero y)
 \\
 \label{d-derivation}
\partial(fg) & = & f\opm \partial(g) + g\opm \partial(f) 
 \\
 v \zero \partial(f)  & = & \partial(\pi(v)(f)) 
 \\
 \label{va-1d}
 v \one \partial(f)   & = & \pi(v)(f)
 \end{eqnarray}
for $v,x,y\in\cA$, $f,g\in\cO_X$.
The map $\pi$ is called the {\em anchor} of $\cA$.

If $\cV = \bigoplus_{n \ge 0} \cV_n$ 
is a (graded) sheaf of vertex algebras with $\cV_0 = \cO_X$, 
then  $\cA = \cV_1$ is a  vertex algebroid with $\de$ equal to the translation operator 
and $\pi$ sending $x\in \cV_1$ to the derivation $f \mapsto x \zero f$.

\subsubsection{Associated Lie algebroid} 
Recall that  a {\em Lie algebroid} 
  is a sheaf of $\cO_X$-modules $\cL$ equipped with a Lie algebra bracket $[, ]$
 and a morphism
$
\pi: \cA \to \cT_X
$
of Lie algebra and $\cO_X$-modules
called {\em anchor}
that satisfies
   $
   [x, ay] = a[x,y] + \pi(x) (a) y
   $,
    $x,y \in \cA$, $a\in \cO_X$.

 \smallskip
 
 If $\cA$ is a vertex algebroid, then 
 the operation $\zero$ descends to that on
 $\cL_{\cA}= \cA / \cO_X \opm \de \cO_X$ 
 and makes it into
  a Lie algebroid, with the anchor induced by that of $\cA$.
  $\cL_{\cA}$  is called {\em the associated Lie algebroid of} $\cA$.

\subsubsection{}
  A vertex (resp., Lie) algebroid is {\em transitive},   if its  anchor map $\pi$ is surjective.

 Being a derivation, see (\ref{d-derivation}),
  $\de: \cO_X \to \cA$ lifts to $\Omega^1_X \to \cA$.
 It follows from (\ref{va-1d}) that if $\cA$ is transitive, then
 $\Omega^1_X \iso \cO_X \opm \de \cO_X$ and $\cA$ fits into an exact sequence
 $$
 0 \To \Omega^1_X \To \cA  \To \cL \To 0,
 $$
$\cL = \cL_{\cA}$ being an extension
  $$
  0 \To \fh(\cL)  \To \cL   \To \cT_X  \To 0
  $$
 where $\fh(\cL) := \ker ( \cL \stackrel{\pi}{\To}  \cT_X )$ is an $\cO_X$-Lie algebra. 
 
\smallskip

 Note that the pairing $\one$ on $\cA$ induces 
 a symmetric   $\cL_{\cA}$-invariant $\cO_X$-bilinear pairing on $\fh(  \cL_{\cA})$
 which will be denoted by $\pair{ , }$.
 
 We regard the pair $(\cL_{\cA}, \pair{, } )$ as "classical data" underlying the vertex algebroid $\cA$.

\subsubsection{Truncation and vertex enveloping algebra functors}

There is an obvious  truncation functor $$t: \cV ert \to \cV
ert_{\leq 1} $$ that assigns to every vertex algebra a 1-truncated
vertex algebra.
  This functor admits a left adjoint \cite{GMS}
  $$ u:  \cV ert_{\leq 1}  \to \cV ert $$
   called a {\em vertex enveloping algebra functor}.

These functors have evident sheaf versions. 
  In particular, one has the functor
\begin{equation}
U: \cV ert\calA lg \To Sh   \cV ert 
\end{equation}
from the category of vertex algebroids to the category of  sheaves of vertex algebras.

\subsection{Chiral differential operators}
 \label{Chiral_differential_operators} 
     A sheaf of vertex algebras $\cD$ over $X$ is called a
     {\em sheaf of chiral differential operators}, CDO for short, 
    if $\cD$ is the vertex envelope of a vertex algebroid 
     $\cA$  that fits into an exact sequence of $\Cplx$-vector spaces
     $$
      0 \To \Omega^1_X \To   \cA   \To   \cT_X   \To   0.
     $$

A sheaf of chiral differential operators  does not  exist over any $X$,
but it does exist locally on any smooth $X$.

 To be more precise, a smooth affine variety $U= \Spec A $ will
  be called {\em suitable for chiralization}
 if $Der(A)$ is a free $A$-module admitting an abelian frame $\{\tau_1,...,\tau_n\}$.
 In this case
  there is a CDO over $U$, which is uniquely determined by the
  condition that
  $(\tau_i)_{(1)}(\tau_j)=(\tau_i)_{(0)}(\tau_j)=0$.
  Denote this CDO by $D^{ch}_{U,\tau}$.

\begin{thm}
\label{class-cdo-local}
 Let $U=\Spec A$ be suitable for
chiralization with a fixed abelian frame $\{\tau_i\}\subset Der A$.

(i) For each closed 3-form $\alpha\in\Omega^{3,cl}_A$ there is a CDO
over $U$ that is uniquely determined by the conditions
\[
(\tau_i)_{(1)}\tau_{j}=0,\;(\tau_i)_{(0)}\tau_{j}=\iota_{\tau_i}\iota_{\tau_j}\alpha.
\]
Denote this CDO by $\cD_{U,\tau}(\alpha)$.

(ii) Each CDO over $U$ is isomorphic to $\cD_{U,\tau}(\alpha)$ for
some $\alpha$.

(iii) $\cD_{U,\tau}(\alpha_1)$ and $\cD_{U,\tau}(\alpha_2)$ are
isomorphic if and only if there is $\beta\in\Omega^{2}_A$ such that
$d\beta=\alpha_1-\alpha_2$. In this case the isomorphism is
determined by the assignment
$\tau_i\mapsto\tau_i+\iota_{\tau_i}\beta$.
\end{thm}

If $A=\pole[x_1,...,x_n]$, one can choose $\partial/\partial x_j$,
$j=1,...,n$, for an abelian frame and check that the beta-gamma
system $M$ of sect. \ref{beta-gamma-system} is a unique up to
isomorphism CDO over $\pole^n$. A passage from $M$ to Theorem
\ref{class-cdo-local} is accomplished by the identifications
$b^j_01=x_j$, $a^j_{-1}1=\partial/\partial x_j$.

\bigskip

\section{Universal  twisted cdo}

In this section we recall the definition of the sheaf $\twcdo$
of twisted chiral differential operators (TCDO) \cite{AChM}
corresponding to a given CDO $\cD^{ch}$
on a smooth projective variety $X$.

\subsection{Twisted differential operators.}
\label{Twisted_differential_operators} A sheaf of twisted
differential operators (TDO) is a sheaf of filtered $\cO_X$-algebras
such that the corresponding graded sheaf is (the push-forward of)
$\cO_{T^*X}$, see [BB2]. The set of isomorphism classes of such
sheaves is in 1-1 correspondence with $H^1(X,\Omega^{[1,2>}_X)$.
Denote by $\cD^\la_X$ a TDO that corresponds to $\la\in
H^1(X,\Omega^{[1,2>}_X)$. If $\text{dim}
H^1(X,\Omega^{[1,2>}_X)<\infty$, then it is easy to construct a
universal TDO, that is to say, a family of sheaves with base
$H^1(X,\Omega^{[1,2>}_X)$ so that the sheaf that corresponds to a
point $\la\in H^1(X,\Omega^{[1,2>}_X)$ is isomorphic to $\cD^\la_X$.
The construction is as follows.

{\em Let $X$ be a  smooth  projective variety}. Then 
$\dim H^1(X,\Omega^{[1,2>})<\infty$
and  there exists an affine cover $\frak{U}$ so that
$\check{H}^1(\frak{U},\Omega^{[1,2>})=H^1(X,\Omega^{[1,2>})$.

Let $\Lambda =\check{H}^1(\frak{U},\Omega^{[1,2>}) $.
We fix a lifting 
$
 \check{H}^1(\frak{U},\Om{1}{2}{X}) \To    \check{Z}^1(\frak{U},\Om{1}{2}{X})
 $
 and identify the former with the subspace of the latter defined by this lifting. Thus, each 
 $ \la \in \Lambda$ is a pair of cochains 
 $\la = (  (\la^1_{ij} ),   ( \la^2_{i}))$ with 
 $
  \la^1_{ij}  \in \Omega^1(U_i \cap U_j),
 $
 $
 \la^2_{i}  \in \Omega^{2, cl}(U_i),
 $
 satisfying
 $
 d_{DR}  \la^{1}_{ij}   =  d_{\check{C}}   \la^{2}_{i}
 $
 and 
 $
  d_{\check{C}}  \la^1_{ij}   = 0.
 $

For $\la = (\la^1_{i},   \la^2_{ij}) \in \Lambda$ denote $\cD^{\la}$ 
the corresponding sheaf of  twisted differential operators. 
One can consider $\cD^{\la}$ as 
an enveloping algebra   
of the Lie algebroid $\cT^{\la}  = \cD^{\la}_1$ \cite{BB2}.
As an $\cO_X$-module, $\cT^{\la}$ is an extension
\nc{\dar}{\downarrow}
\nc{\One}{{\mathbf{1}}}
$$
0 \To    \cO_X \One     \To    \cT^{\la}     \To \cT_X     \To  0
$$
given by $(\la^1_{ij})$. 
The Lie algebra structure on $\cT^{\la}_{U_i}$ is given by 
$
[\xi, \eta]_{\cT^\la}  =  [\xi, \eta]  +  \ixi \ieta \la^2_{i}
$ and 
$[\One, \cT^{\la}_{U_i}]= 0$.

\smallskip

 Let $\set{  \la^*_i  }$ and $\set{  \la_i  }$ are dual bases 
of $\Lambda^* = H^1(X, \Om{1}{2}{X})$ and $\Lambda$ respectively.
Denote by $k$  the dimension of $\Lambda$.


Define $\cT^{tw}$  to be  an abelian extension
$$
0 \To \cO_X \tensor \Lambda^*     \To   \cT^{tw}_X  \stackrel{\pi}{\To} \cT_X   \To 0
$$
such that $[\Lambda^*, \cT^{tw} ] =0$ and there exist
sections  
$\nabla_i: \cT_{U_i} \to \cT^{tw}_{U_i}$ of $\pi$
satisfying
\begin{gather}
\label{univTDO-A}
\nabla_j(\xi) - \nabla_i(\xi) =  \sum_{r} \iota_{\xi} \la^{(1)}_r(U_{ij}) \la^*_r
\\
\label{univTDO-curvature}
[\nabla_i(\xi), \nabla_i(\eta)]- \nabla_i ([\xi,\eta]) =  \sum_r   \iota_{\xi}  \iota_{\eta} \la^{(2)}_r (U_{i})  \la^*_r
\end{gather}

It is clear that the pair 
$
(\cT^{tw}, \cO_X \tensor \Lambda^* \hookrightarrow \cT^{tw})
$ 
 is independent of the choices made.

	We call the universal enveloping algebra 
	$\cD^{tw} = U(\cT^{tw})$
	\textit{the 	universal sheaf of twisted differential operators}.

	\subsection{A universal twisted CDO}
	\label{A_universal_twisted_CDO}

	Let $ch_2(X)=0$ 
	It is proved in \cite{GMS} that $X$ carries at least one CDO.
	To each CDO $\cD^{ch}_X$ one
	can attach a {\em universal twisted CDO}, $\cD^{ch, tw}_X$,
	a sheaf of vertex algebras whose "underlying" Lie algebroid is
	$\cD^{tw}_X$.
	 Let us place ourselves in the situation of the previous section, where we had a fixed affine
	cover ${\frak U}=\{U_i\}$ of a projective algebraic manifold $X$,
	 dual bases $\{\lambda_i\}\in
	H^1(X,\Omega^{[1,2>}_X)$, 
	$\{\lambda_i^*\}\in
	H^1(X,\Omega^{[1,2>}_X)^*$, and a lifting
	$H^1(X,\Omega^{[1,2>}_X)\rightarrow Z^1({\frak
	U},\Omega^{[1,2>}_X)$.

        We can assume that $U_i$ are suitable for chiralization. 
         Let us fix, for each $i$, an abelian basis
	$\tau^{(i)}_1,\tau^{(i)}_2,...$ of $\Gamma(U_i,\cT_X)$.
	 Then the CDO $\cD^{ch}$ is given by a
	collection of 3-forms $\alpha^{(i)}\in\Gamma(U_i,\Omega^{3,cl}_X)$
	(cf. sect.\ \ref{Chiral_differential_operators},
	Theorem~\ref{class-cdo-local}) 
	and transition maps
	$
	  g_{ij} : \cD^{ch}_{U_j}|_{U_i \cap U_j } \to  \cD^{ch}_{U_i}|_{U_i \cap U_j }.
	$
	Let us as well fix splittings $\cT_{U_i } \hookrightarrow \cD^{ch}_{U_i}$
	and view $g_{ij}$ 
	as maps 
	$
	 g_{ij}: (\cT_{U_j } \oplus \Omega^1_{U_j} )|_{U_i \cap U_j } \to  (\cT_{U_i } \oplus \Omega^1_{U_i})|_{U_i \cap U_j }
	 $	
	
	The   {\em universal sheaf of twisted chiral differential operators}
	 $\cD^{ch,tw}_X$  (TCDO for short) corresponding to $\cD^{ch}_X$
	is 
	a vertex envelope of the $\cO_X$-vertex algebroid $\cA^{tw}$ determined by the following:
	\begin{itemize}
	\item
	   the associated Lie algebroid of $\cA^{tw}$ is $\cT^{tw}_X$;
	\item
	  there are embeddings  $\cT^{tw}_{U_i} \hookrightarrow \cA_{U_i}$, 
	   such that
	  \begin{gather}
	  \nonumber
	  	   \tau^{(i)}_l \zero \tau^{(i)}_m   
	      =   \iota_{  \tau^{(i)}_l }  \iota_{  \tau^{(i)}_m } \alpha^{ (i) }  
	       + \sum  \iota_{  \tau^{(i)}_l }  \iota_{  \tau^{(i)}_m }   \lambda^{(2)}_k (U_i) \lambda^*_k    
	   \\   
	    \nonumber
	       \lambda^*_k \one x =0, \ x\in \cA^{tw}. 
	  \end{gather}
	 \item 
	   the transition function  from $U_j$  to $U_i$ is given by
	   \begin{equation}
	     \label{twcdo_trans_fun}
	   g_{ij}^{tw} (\xi)  =  g_{ij}(\xi)  - \sum \iota_{\xi} \lambda_k^{(1)} (U_i \cap U_j) \lambda^*_k 
	   \end{equation}
	\end{itemize}
	See \cite{AChM} for a  detailed construction.

   \bigskip

\subsection{Locally trivial twisted CDO}   
   
Observe that there is an embedding
\begin{equation}
\label{emb_coho}
 H^{1}(X,\Omega^{1,cl}_X)\hookrightarrow H^{1}(X,\Omega^{[1,2>}_X)
\end{equation}
The space $H^{1}(X,\Omega^{1,cl}_X)$ classifies {\em locally
trivial} twisted differential operators, those that are locally
isomorphic to $\cD_X$. Thus for each $\lambda\in
H^{1}(X,\Omega^{1,cl}_X)$, there is a unique up to isomorphism TDO
$\stackrel{\circ}{\cD}^{\lambda}_X$ such that for each sufficiently
small open $U\subset X$, $\stackrel{\circ}{\cD}^{\lambda}_X|_U$ is
isomorphic to $\cD_U$. Let us  see what this means at the level of
the universal TDO.

In terms of Cech cocycles the image of embedding (\ref{emb_coho}) is
described by those $(\lambda^{(1)},\lambda^{(2)})$, see section
\ref{Twisted_differential_operators}, where $\lambda^{(2)}=0$, and this
forces $\lambda^{(1)}$ to be closed.  Picking a collection of such
cocycles that represent a basis of $H^{1}(X,\Omega^{1,cl}_X)$ we can
repeat the constructions of sections
\ref{Twisted_differential_operators} and
\ref{A_universal_twisted_CDO} to obtain
a sheaf 
$\stackrel{\circ}{\cD}^{ch,tw}_X$.
It is defined
by gluing pieces isomorphic (as vertex algebras) to
$\cD^{ch}_{U_i} \otimes H_X$ with transition functions as in
(\ref{twcdo_trans_fun}); here $H_X$ is the vertex algebra of differential
polynomials on
 $H^{1}(X, \Omega^{1,cl}_X)$. 
 We will call the 
 sheaf
   $\stackrel{\circ}{\cD}^{ch,tw}_X$ the {\em universal locally trivial 
   sheaf of twisted chiral differential
operators.}

\subsubsection{Example}
\label{Example:_prline.} Let us construct  a sheaf of TCDO on
 $X=\prline$. We have
$\prline=\pole_{0}\cup\pole_{\infty}$, a cover ${\frak
U}=\{\pole_{0},\pole_{\infty}\}$, where $\pole_0$ is $\pole$ with
coordinate $x$, $\pole_{\infty}$ is $\pole$ with coordinate $y$,
with the transition function $x\mapsto 1/y$ over
$\pole^*=\pole_{0}\cap\pole_{\infty}$.

Defined over $\pole_{0}$ and $\pole_{\infty}$ are the standard CDOs,
$\cD^{ch}_{\pole_0}$ and $\cD^{ch}_{\pole_{\infty}}$. The spaces of
global sections of these sheaves are  polynomials in
$\partial^{n}(x)$, $\partial^{n}(\partial_x)$ (or $\partial^{n}(y)$,
$\partial^{n}(\partial_y)$ in the latter case), where $\partial$ is
the translation operator, so that, cf. sect.\
\ref{Chiral_differential_operators},
$$
(\partial_x)_{(0)}x=(\partial_y)_{(0)}y=1.
$$
There is a unique up to isomorphism CDO on $\prline$,
$\cD^{ch}_{\prline}$; it is defined by gluing $\cD^{ch}_{\pole_0}$
and $\cD^{ch}_{\pole_{\infty}}$ over $\pole^*$ as follows
\cite{MSV}:
\begin{equation}
 \label{gluing_untwisted_p1}
x\mapsto 1/y,\; \partial_x\mapsto
(-\partial_{y})_{(-1)}(y^2)-2\partial(x).
\end{equation}

The  twisted version is as follows.

Since $\dim \prline=1$,
$$
H^1(\prline,\Omega^{1}_{\prline}\rightarrow\Omega^{2,cl}_{\prline})=
H^1(\prline,\Omega^{1, cl}_{\prline})
$$
 so all twisted CDO on $\prline$ are locally trivial. 
Furthermore,
$H^1(\prline,\Omega^{1,cl}_{\prline})=\pole$ and is spanned by the cocycle
$\pole_0\cap\pole_{\infty} \mapsto dx/x$. 
  We have
  $H_{\prline}=\pole[\lambda^*,\partial(\lambda^*),....]$. 
Let
$\cD^{ch,tw}_{\pole_0}=\cD^{ch}_{\pole_0}\otimes H_{\prline}$,
$\cD^{ch,tw}_{\pole_{\infty}}=\cD^{ch}_{\pole_{\infty}}\otimes H_{\prline}$ and define $
\cD^{ch,tw}_{\prline}$ by gluing $\cD^{ch,tw}_{\pole_0}$ onto $ \cD^{ch,tw}_{\pole_{\infty}}$
via
\begin{equation}
 \label{gluing_twisted_p1}
\lambda^*\mapsto\lambda^*,\;x\mapsto 1/y,\;\partial_x\mapsto
-(\partial_y)_{(-1)}y^2-2\partial(y) +y_{(-1)}\lambda^*.
	\end{equation}

 All of the above is easily verified by direct computations, cf. \cite{MSV}.

\section{Main  result}
\label{Main result - section}
In this section we prove the existence of a filtration for modules $\cM$ 
 that satisfy certain integrability condition
 and use it to prove a generalization of Theorem 5.2. in \cite{AChM}.
 
We work in a setup that is slightly more general than that of TCDO
 with an intention  to apply
 the results to
certain deformations of TCDO defined in \cite{Ch}.

Let $\cV$ be a sheaf of vertex algebras. 
	We will call a sheaf of vector spaces $\cM$
	a $\cV$-{\em module} if
    for each open $U \subset X$ 
	$\cM (U)$
	is a
	$\cV (U)$-module and
    the restriction morphisms
	$\cM(U)\rightarrow   \cM(U' )$, $U' \subset U$, 
	 are  $\cV(U)$-module morphisms, with  the
	 $\cV(U)$-module structure
	  on $\cM(U')$  given by  pull-back.

\smallskip

Recall that to each graded vertex algebra $V$ one can associate 
$\Lie V$,  the Borcherds Lie algebra of $V$ (cf. section \ref{BorcherdsLie}).

We say that a $\cV$-module $\cM$ is {\em  half-integrable}
 iff each point $x\in X$ has a neighbourhood $U$ 
 with $\cM(U)$  a $\Lie(\cV(U))_+$-integrable module, that is to say, 
 $\Lie(V)_+$ acts on $M$ by locally nilpotent operators.

 We call the {\em center}  of $\cV$  the subsheaf $\cZ (\cV) \subset \cV$   with stalks $\cZ_{X,x} = Z(\cV_{X,x})$, $x\in X$
 where $Z(V)$ denotes the center of $V$.

\bigskip
 
Define the category $\halfint(\cV)$ as the full subcategory of the category of $\cV$-modules
$\cM$ such that
\begin{enumerate}
\item   
      $\cM$ is  half-integrable  
\item 
$
h_n  m =0, \textrm{  for all }  h \in \cZ(\cV), m \in \cM.
$
\end{enumerate}

\smallskip

Our goal in this section is the proof of the Theorem below.
We  work in the analytic topology.
In section \ref{the-algebraic-case-Section}
we present a version of this result which works
in Zariski topology. 
\begin{thm}
\label{main_result_gen}

 Suppose $\cV$ is a vertex envelope of a
 transitive vertex algebroid $\cA$ 
 whose associated Lie algebroid $\cL = \cA / \Omega^1_X$
  is locally free of finite rank 
 and
 fits into an exact sequence
\begin{equation}
\label{Lie_algd_ses2}
0 \To \fh \To  \cL \To \cT_X \To 0
\end{equation}
in which $\fh$ is an {\em abelian} $\cO_X$-Lie algebra.
Let  $\cM \in \halfint(\cV)$. Then 

(1) \ $\cM$ is generated by the subsheaf
 $$
  \Sing \cM = \set{ m \in \cM : \ v_n m = 0  \textrm{ for all } v\in \cV,  \,  n>0 }
 $$

(2) \ 
There is a filtration
 $\set{\cM_i}_{i\ge 0}$
 on $\cM$ with $\cM_0 = \Sing \cM$,
compatible with the grading of $\cV$.
\end{thm}

\subsection{Proof}

\subsubsection{The center}
First of all, let us describe the center of any such $\cV$.
For that, let us look more closely at  (\ref{Lie_algd_ses2}). Since $\fh$ is abelian, it possesses
        the structure of a $\cD_X$-module. It is
	  locally free as an  $\cO_X$-module 
	and therefore, is of the form 
	$
	\fh = \cO_X \tensor_\Cplx   \fh^{\nabla} 
	$
	where $ \fh^{\nabla}$ denotes the subsheaf of horizontal sections of $\fh$.

 The pairing $\one$   induces an $\cO_X$-bilinear $\cL$-invariant  symmetric
  pairing $\pair{,}: \fh \times \fh \to \cO_X$ which restricts to
$
\pair{ , }: \ \fhn \times \fhn \to \Cplx_X
$.

\smallskip

Let
$
\cz = \set{  h \in \fhn \, : \ \pair{ h, \fhn} = 0}
$.
 Thus $\cz$ is the kernel of $\pair{,}$ restricted to $\fhn$.

  \smallskip

  \begin{lem}
  \label{lem-central-lifting}
  There exists a lifting  $s: \cz \to \cA$, i.e. $\pi s = \id_{\cz}$, 
   such that $s(h)$, $h \in \cz$ 
   generate a  subalgebra in $\Gamma(X, \cV)$
   that is central in every $\Gamma(U,\cV)$, $U \subset X$.
   Moreover, such a lifting is unique.
   \end{lem}
   \smallskip
   
   {\em Proof.}
   Let $s': \cz \to \cA$ be any (local)  lifting. 
   Fix a basis ${\tau_i}$ in $\cT(U)$ and the dual basis $\set{\omega_i}$ in $\Omega^1(U)$.
    Fix arbitrary lifts $\ttau_i$ of $\tau_i$ in $\cA$.
    For $s'(h)$ to be central, $s'(h) \one \ttau_i$ must be zero  for all $i$.
    This may fail, and so we are forced to change $s'$. It is clear that
    $$
    s(h) = s'(h) - ( \ttau_i \one s'(h))  \omega_i, \ \ h\in \cz
    $$ 
    satisfies the desired condition
    $s(h) \one \ttau_i  =0$
    for all $h\in \cz^*$ and all basis elements $\tau_i$; 
    furthermore, the latter condition determines $s(h)$ uniquely.
    Therefore,  $s$ is unique and extends globally. 
    
    It remains to show that 
    $\ttau_i \zero s(h) =0$ for all $i$
   
    Now,  $[\tau_i, h]=0$ in $\cL$ (since $\cz \subset \fhn$).
     Therefore
    $\ttau_i \zero s(h)$ is in $\Omega^1$.  
    However, for any $1\le j \le \dim X$
    $$
   ( \ttau_i \zero s(h) ) \one \ttau_j  =  \ttau_i \zero (s(h)  \one \ttau_j ) - s(h) \one (\ttau_i  \zero \ttau_j) = 0
    $$
    Thus   $\ttau_i \zero s(h)$ must be zero.      
    $\qed$

\bigskip

  	Let
	 $U \subset X$ be suitable 
	 for chiralization (cf. \ref{Chiral_differential_operators})
   Choose an abelian basis $\set{\tau_i}$ of $\cT (U)$,
   the dual basis $\set{ \omega_i}$ in $\Omega^1(U)$
   and  a  basis $\set{h_k}$ of   $\fhn(U)$. 

 Choose any lifting $\cL \to \cA$ extending that of Lemma \ref{lem-central-lifting}
 and identify $\tau_i$ and $h_k$ with the corresponding lifts.

Let $V$  be the vertex enveloping algebra of  $\cA(U)$. 
$V$ is generated by $\cO_{U}$,  vector fields $\{ \tau_i \}$    
and elements $h_k$.

The following relations hold in $V$.
\begin{eqnarray}
\label{first_rel}
[\omega_{i,n}, \omega_{j,l}]  &  = & 0 \\
\label{second_rel}
[\tau_{i,n}, \omega_{j,l}]  &= & n \delta_{i,j}  \delta_{n+l, 0}  \id  \\
\label{third_rel}
\left[\tau^{}_{i,n}, \tau_{j,l} \right]  & = &  \alpha_{n+l}  + \sum_k    \left( f_k \opm \lambda^*_k \right)_{n+l}    \\
\label{fourth_rel}
\left[h_{k,n} , h_{l,m}   \right] & = & 
 n \pair{h_k, h_l} \delta_{n+m, 0}\id
\\
\label{fifth_rel}
[\tau_{i, n}, h_{r,m} ] & = &  \beta_{n+m} 
\end{eqnarray} 
where $\alpha \in \Omega^1(U)$, 
    $f_k\in \cO(U)$  and $\beta\in  \Omega^1(U)$
  may  depend on $\tau_i$, $\tau_j$ and $h_{r }$ resp.

 Using (\ref{second_rel} - \ref{fifth_rel}), 
 it is easy to show that the subalgebra generated by $\Gamma(U, s(\cz))$,
 see Lemma \ref{lem-central-lifting},
 is all of $Z(V)$.

\subsubsection{} 
Let us now describe the strategy of proving Theorem \ref{main_result_gen}, (1). 
 The statement is local and we continue working on a suitable 
	 for chiralization subset $U$.

Let $M = \cM(U)$.

Denote
$\Sing M = 
  \set{ m \in M : v_n m = 0 \textrm{ for all  } v\in V , n>0 }$

Introduce the filtration:
\begin{equation}
\label{filtrationM}
 0 \subset \Moh \subset \oinv  \subset M
\end{equation}
where 
$$
\oinv =  \set{ m\in M:  \omega_{i, n}m = 0  \ \textrm{ for all }  n>0, i  }
$$
$$
\Moh =  \set{ m\in \oinv:  h_{r,n}m = 0 \textrm{ for all }  n>0,  r }
$$

Let $\byzero$ be the submodule of $M$ generated by $\Sing M$.  

We show, step by step, that each of the terms in the filtration (\ref{filtrationM})
is generated by  $\Sing M$.

\smallskip

\smallskip

\subsubsection{Step 1}
We show that $\byzero$ contains the subspace $\Moh$.

\smallskip

Let us introduce a filtration on $\Moh$ indexed by 
 functions  
 $$
 \dd : \   \set{1,2,\dots, \dim X}\times \ZZ_+  \to \ZZ_+
 $$
vanishing at all but finitely many pairs $(i, n)$. 

Call any such function a {\em degree vector}.

 For a degree vector $\dd$  define 
 \begin{equation}
 \Moh  (\dd ) =  \bigcap_{i, n>0}   \ker  \tau_{i,n}^{d(i,n)+1} .
 \end{equation}

It is clear from definitions that
 $
 \Moh = \bigcup \Moh (\dd)
 $ 
 and 
{
 $\Moh(0) = \Sing M$. 
}

Let us show that each
 $\Moh (\dd) $  
 is a subspace of   
 $\byzero$
using the induction on the {\em length} 
$\abs{\dd} = \sum_{i,n} d(i, n)$ of $\dd$. 
The base of induction 
 is established in the line above.

  Let $\dd  \ne 0 $   
  be a degree vector. 
  Suppose 
 $ \Moh (\dd') 
 \subset  \byzero$
 for all 
 $\dd'$ 
 of smaller length.
 
 \smallskip
 
 Fix $(i,n)$ such that $k: = d(i,n)  > 0$
 and let
 $\dd'$
 be equal to $\dd$ everywhere  except at $(i,n)$ where 
 $\dd'(i,n)  = \dd(i,n)  - 1$. 
By induction assumption, 
$\Moh (\dd') \subset \byzero$.

 \smallskip
 
 Let
 $
m \in \Moh(\dd) 
$.
 Since $m \in \ker \tau_{i, n}^{k+1}$,
 one has  
 \begin{equation}
 \label{omega_i_tau}
  0 =   \omega_{i, -n} \tau_{i, n}^{ k + 1} m    =  
      \tau_{i, n}^{ k  } ( - nk m  +  \omega_{i, -n} \tau_{i, n} m ) 
 \end{equation}

  Introduce the elements
\begin{equation}
 \begin{split}
  m' &= \tau_{i, n}m, \\
m''  &= nk m  - \omega_{i, -n} m'.
\end{split}
\end{equation}
Then 
\begin{equation}
\label{m_stepone}
m = \frac{1}{nk} (m'' + \omega_{i, -n} m')
\end{equation}
Hence, to show $m \in \byzero$ 
it suffices to show that $m'$, $m''$ lie in $\byzero$.

 \smallskip

Let us  show 
 $m' \in \Moh( \dd')$.

One has $\tau_{i,n}^{k} m'  = \tau_{i,n}^{k+1}m = 0$.  To show 
$
\tau_{j,l}^{\dd(j,l) +1} \tau_{i,n} m = 0 
$
it  suffices to show
$$
[ \tau_{j,l},  \tau_{i,n}]  \tau_{j, l}^{q} m =0, \ \ q\ge 0.
$$
But that follows from 
the  fact that $\tau_{j,l} \Moh \subset \Moh$, which is a consequence of (\ref{second_rel}), (\ref{fifth_rel}). 

\smallskip

Finally,  to see
 that  $m''$   is in $F^{\dd ' }$,
we need to check
that $\tau_{j, l}^{\dd'(j, l)  + 1 } m'  =0$ for all $(j,l)$.
For $(j, l) \neq (i, n)$ 
this follows immediately from (\ref{second_rel}) and 
the fact that $m' \in F^{\dd'}$. The case $(j,l)=(i,n)$ is clear due to  (\ref{omega_i_tau}).

\bigskip

\begin{rem}
Note that (\ref{omega_i_tau})
is a particular case of the following observation used in Steps 2 and 3 as well and originating in Kashiwara's lemma.

{\em 
\label{GKL_lem}
Let $A$ and $B$ be linear operators on a space $V$ such that $[A,B]$ commutes with $A$. 

Suppose $A^{n+1} m=0$ for $m \in V$, $n\ge 0$.
Then 
\begin{equation}
\label{GKL_eq}
 A^n \cdot (n [B,A] + BA )m  =0
\end{equation}
  }

Indeed,
$
0 = BA^{n+1} m = [B, A^n] Am  + A^n \, BAm =  A^n (n [B,A] + BA) m \  \qed
$

\end{rem}

\bigskip

Steps 2 and 3 are proved in essentially the same way. 
The reader may safely skip the rest of the proof.
However,  for the sake of completeness we will keep the same level of detail.

\smallskip

\subsubsection{Step 2}      
 We show
{ $\oinv \subset \byzero$.}

In the case  $\pair{ , } = 0$
 the assumptions  of the theorem imply 
 $h_{r,n} =0$ on $M$  so that $\oinv = \Moh$  and there's nothing to prove.

Let us prove the claim
 in the case  
{  $\cz $  is a  proper subset of     $\fhn$. }

 Let
 $\set{b_r}$ 
 be an orthonormal basis for some complement 
 $\fhn_{nondeg}$ 
 to
  $\cz$ in  $\fhn$.
 In particular
 \begin{equation}
 \label{fourth_rel_orth}
\left[b_{k,n} , b_{l,m}   \right]  =   n  \delta_{k, l} \delta_{n+m, 0} \id 
 \end{equation}
 
\smallskip

Introduce the filtration 
$$
\oinv(\dd) = \set{ m \in \oinv \ : \  (b_r)_n^{\dd(r, n ) + 1} m = 0, \forall n>0 \ \forall r }
$$
where
$
\dd:   \{1.. \dim \fhn_{nondeg}   \} \times \ZZ_+  \to \ZZ_+
$
is a degree vector.
\smallskip

Then, clearly, 
$\oinv(0) = \Moh$
and 
  $\oinv  = \bigcup_{\dd} \oinv(\dd)$.

  Suppose that
   $\dd  \neq {\mathbf{0}}$ is a degree vector
and suppose 
 $\oinv(\dd')$ is a subspace of 
  $\byzero$
 {\em for all }
 $\dd'$ 
 of smaller length. 
 Let us show $\oinv (\dd) \subset \byzero$

Fix some $i, n$ such that $d(i,n)>0$ and define $\dd'$ by 
 $\dd' (j,l) = \dd(j, l) - \delta_{ij} \delta_{nl}$.
 Let $k=  \dd(i,n)$.
 
 \smallskip
 
{ Let $m \in \oinv (\dd)$. }
 
 We have (cf. (\ref{GKL_eq}))  
 \begin{equation}
 \label{kash_lambda}
 0    =  (b_{i, n})^{k}  \left( -nk  m + b_{i, -n} b_{i, n}  m \right) 
 \end{equation}

  Introduce the elements
\begin{equation}
 \begin{split}
  m' &= b_{i, n } m, \\
m''  &= -nk   m + b_{i, -n} m'    
\end{split}
\end{equation}

Since $\omega$'s commute with $b_r$'s,
 these elements are in $\oinv$ whenever $m$ is.
We wish to show that they are in fact in $\oinv(\dd')$ and therefore, in $\byzero$.  
This would imply  
that 
\begin{equation}
 \label{m_steptwo}
m = \frac{1}{nk} ( - m'' + b_{i, -n} m')
\end{equation}
is in $\byzero$ as well.

\smallskip

It is clear from (\ref{fourth_rel_orth}) that  $m' \in \oinv(\dd')$. 

Let us show that  $m'' \in   \oinv(\dd')$. 
We need to show 
$
m'' \in \ker (b_{s,m}  )^{\dd'(s,m) +1}
$
for all $s$ and all  $m>0$; but
 this follows from (\ref{fourth_rel_orth}) in case $(s, m) \ne (i,n)$ and from 
(\ref{kash_lambda})
in case $(s,m) = (i, n)$. Thus  $m'' \in   \oinv(\dd')$.

\bigskip

\subsubsection{Step 3} 
We complete the proof by showing $M \subset \byzero$.

Introduce a filtration on $M$ indexed by degree vectors
   $\dd : \   \set{1,2,\dots, \dim X}\times \ZZ_+  \to \ZZ_+$
where for each $\dd$ we define
 \begin{equation}
 M(\dd ) =  \bigcap_{i, n>0}   \ker ( \omega_{i,n}^{d(i,n)+1}).
 \end{equation}

Clearly, $M = \bigcup M(\dd)$ and 
{
 $M(0) = \oinv$. 
}

In order to show that $M(0) \subset \byzero$ implies
 $M(\dd) \subset \byzero$  for all $\dd$
   we will use induction on the length $\abs{\dd} = \sum_{i,n} d(i, n)$.
  Thus,  we suppose that
   $\dd  \neq {\mathbf{0}}$ is a degree vector
and 
 $M(\dd')$ is a subspace of   $\byzero$
  for all 
 $\dd'$ 
 of smaller length. 
 Let us show $M(\dd) \subset \byzero$

Fix some $i, n$ such that $d(i,n)>0$ and define $\dd'$ by 
 $\dd' (j,l) = \dd(j, l) - \delta_{ij} \delta_{nl}$.
  Let $k=  \dd(i,n)$.

 \smallskip
 
   Let $m \in M(\dd)$. 
 
 Since $m \in \ker \omega_{i, n}^{ k + 1}$, one has (cf.  (\ref{GKL_eq}) )
 \begin{equation}
 \label{ tau_i_n_omega_i_n}
  0 = \tau_{i, -n} \omega_{i, n}^{k+ 1} m =   \omega_{i, n}^{  k  } (- nk m  + \tau_{i, -n} \omega_{i, n} m ) 
 \end{equation}

It follows from
   (\ref{first_rel} -\ref{third_rel}) 
   and
    (\ref{ tau_i_n_omega_i_n})
   that the elements
\begin{equation}
 \begin{split}
  m' &= \omega_{i, n}m, \\
m''  &= nk m  - \tau_{i, -n} m'
\end{split}
\end{equation}
belong to $M(\dd')\subset \byzero$.
Therefore, 
\begin{equation}
 \label{m_stepthree}
	m = \frac{1}{nk} (m'' + \tau_{i, -n} m')
\end{equation}
is also in $\byzero$.  $\qed$

\bigskip

 \subsubsection{Proof of (2)} 
Let us work locally as in the proof of part (1)
  and keep all the notation used there.  
  Thus we have
$U\subset X$, $M = \cM(U)$, $V$ is the vertex envelope of $\cA(U)$.

	Define
	  $\Lie(V)_+$ to be the subalgebra of the Borcherds Lie algebra 
	  $\Lie(V)$ spanned by $v\tensor t^n$ with positive $n$;
	  and 
	  let $\Lie(V)_{\leq 0}$ be the image of $V\tensor \Cplx[t^{-1}]$ in $\Lie(V)$.
	Then $\Lie(V) = \Lie(V)_+ \oplus Lie(V)_{\leq 0}$

The enveloping algebras $U(\Lie(V))$,
$U(\Lie(V)_+)$ 
have a natural grading defined by
$\deg v_{n_1}\dots v_{n_k} = n_1  + \dots + n_k$.

Define the following subspaces of $M$:
\begin{equation}
\label{def-filtration-local}
 M_n = \set{ m:\   U(\Lie (V)_+ )_k m = 0, \, \text{ \rm for all } k > n  }
 \end{equation}

It is clear that 
$M_n \subset M_{n+1} $,
 $M_{-1}= \{ 0 \}$ and $M_0 = \Sing M$.

The decomposition 
$
U(\Lie(V)) \iso U(\Lie(V)_{\leq}) \tensor U(\Lie(V)_+)
$
implies 
 that $M_n = \tilde{M}_n$ where
$$
\tilde{M}_n  
= \set{ m\in M : \ 
U(\Lie(V))_k m = 0 \text{ \rm for all } k>n 
}
$$
Note  that  $\tilde{M}_n $ 
is compatible with the grading of $V$, i.e.
$
v_k  \tilde{M}_n \subset \tilde{M}_{n-k}.
$
Hence, the submodule $\byzero$ generated by 
$\Sing M = M_0$ 
is a subset of  the union of   $\tilde{M}_n$. 
 But $\byzero = M$ by part (1) and thus, $M_n$ is an (exhaustive)  filtration.

Since $\cV$ is a graded sheaf, the formula  (\ref{def-filtration-local}) makes sense globally, giving the desired filtration $\cM_n$.  $\qed$

\bigskip

\subsubsection{The algebraic case}
\label{the-algebraic-case-Section}

Theorem \ref{main_result_gen} remains true in Zariski topology, 
as long as one imposes one additional restriction.
Since algebraic $\cD$-modules  
may have non-algebraic solutions,
$\fhn$ may fail to have "right size".
Thus, one has to demand that $\fh$ be equal to $\cO_X \tensor \fhn$.
To be more precise, one has:

\begin{thm}
\label{main_result_alg_case}

 Suppose $\cV$ is a vertex envelope of a
 transitive vertex algebroid $\cA$ 
 whose associated Lie algebroid $\cL = \cA / \Omega$
  is locally free of finite rank 
 and
 fits into an exact sequence
\begin{equation}
\label{Lie_algd_ses2}
0 \To \fh \To  \cL \To \cT_X \To 0
\end{equation}
in which $\fh = \cO_X \tensor \fhn$
for a locally constant sheaf $\fhn$ such that 
$[\cL, \fhn] = 0$.

Let  $\cM \in \halfint(\cV)$. Then:

(1) \ $\cM$ is generated by the subsheaf
 $$
  \Sing \cM = \set{ m \in \cM : \ v_n m = 0  \textrm{ for all } v\in \cV,  \,  n>0 }
 $$

(2) \ 
There is a filtration
 $\set{\cM_i}_{i\ge 0}$
 on $\cM$ with $\cM_0 = \Sing \cM$,
compatible with the grading of $\cV$.
\end{thm}

The proof is identical to that of Theorem \ref{main_result_gen}.

\smallskip

%

%
%
%
%
%
%



\subsection{Equivalence}

Now we prove  a generalization of  Theorem 5.2 in \cite{AChM}
 in the setup of the previous section.
 This is a result  establishing the equivalence between  certain subcategories of $\halfint(\cV)$
 and categories of modules over sheaves of tdo,  based on a version of Zhu correspondence.
 The result is true in both analytic and Zariski topology.

  Let $\cV$ be as in Theorem \ref{main_result_gen} and assume further that
  $\fh = \cO_X  \tensor \fhn$ with $\fhn$ a constant sheaf. 
  By abuse of language we denote its fiber by $\fhn$ and write $(\fhn)^*$
  and $\cz^*$ for vector space duals of 
  $\fhn$ and $\cz := \{ h\in \fhn \, : \ \pair{ h,s } = 0 , s\in \fhn \}$

   Let us fix 
 $$
 \theta \in (\fhn)^* , \ \ \chi = \sum_{n \in \ZZ}  \chi_n z^{-n-1} \in \cz^*((z)).
 $$
  
 Define the category 
 $
\halfintchitheta(\cV)
 $
 to be the full subcategory of $\halfint(\cV)$
 consisting of modules  $\cM$ satisfying:
 
 (1)  \ 
  for all $h \in \fhn$,  $m \in \cM$ 
 \begin{equation}
  h_0 m = \theta(h)  \cdot m, 
 \end{equation}
 
 (2) \ 
 for all $c \in  \cz$, $n \in \ZZ$
 \begin{equation}
      c_{n}  m = \chi_n(c) \cdot m,   
 \end{equation} 

For $\halfintchitheta(\cV)$ to be nonzero, 
an evident compatibility condition has to be satisfied:
\begin{equation}
\label{comp-condition}
\theta |_{\cz} = \chi_0
\end{equation}
Furthermore,
the half-integrability condition dictates
\begin{equation}
\label{comp-condition2}
\chi_n = 0,   \ \ n>0
\end{equation}
for any nonzero module in  $\halfintchitheta(\cV)$. 
In other words, 
 $\chi(z)$ has {\em regular singularity} at 0.
 \smallskip

 Let $\cM \in  \halfintchitheta(\cV)$ be arbitrary.
 By Theorem \ref{main_result_gen}, $\cM$ is a $\ZZ_{\ge 0}$-filtered module
 over a sheaf of (graded) vertex algebras $\cV$. 
 Thus, $\cM_0 = \Sing \cM$ is equipped with an action of  the Zhu algebra of $\cV$, or, rather,
 a sheafified version of it. 
 Furthermore, 
 the action of
 $
 \cZ hu(\cV)
 $
 on
 $
 \cM_0 = \Sing \cM 
 $ 
 factors through
  $
 \cZ hu(\cV)  /  I_\theta 
 $
 where $I_{\theta}$ is the ideal generated by 
 $h - \theta(h)$, $h \in \fhn$.
 As 
 $
 \cZ hu(\cV)   =  U(\cL) 
 $
  ( \cite{AChM}, Theorem 3.1),
  the quotient algebra  
   $
 \cZ hu(\cV)  /  I_\theta 
 $
 is a TDO, to be denoted $\cD^{ \bar{\theta}}_X$.

 \smallskip

 Thus  we have  a functor
\begin{equation}
\label{Sing-scalar}
\Sing : \   \halfintchitheta(\cV)  \to \cM(  \cD^{\bar{\theta} }_X )
\end{equation}

\begin{thm}
\label{equivalence-theorem}
Assume that  
the conditions (\ref{comp-condition}, \ref{comp-condition2}) are satisfied.
Then the functor
(\ref{Sing-scalar})
is an equivalence of categories.
\end{thm}

This theorem is a generalization of Theorem 5.2., \cite{AChM}
and has an almost identical proof, after some modifications;
therefore, we  only list the differences; 
an interested reader may easily supply all the details
following the proof in  \cite{AChM}.
 Recall that the strategy of proof was 
to construct  a left adjoint $\cZ hu^*_{\chi(z)}$
to  (\ref{Sing-scalar})
and to prove that it is a quasi-inverse. 
Our proof proceeds in the same way except for the following:

\begin{itemize}
\renewcommand{\labelitemi}{$-$}
\item
 in the construction of the left adjoint functor
one has to change $H^1 (X, \Omega^{[1,2>}_X)$ to $\cz$

\item
 the   quasi-inverse property 
(Lemma 5.3, {\em op.cit.}) is changed as follows:
one  adds
to the polynomial algebra $P$ 
 indeterminates 
 $L^{r}_{-n}$
  corresponding to 
  $l_{r, -n}$, 
where $\{  l_r \}$ is an orthogonal basis of a complement to $\cz$ in $\fhn$. 
The proof of this new statement  is identical to that of Lemma 5.3.
\end{itemize}
In the rest, one proceeds literally in the same manner as in {\em op.cit.}
$\qed$

\bigskip

As an example, consider the case of TCDO, $\cV =  \twcdo$.
In this case $\cL = \cD^{tw}$, $\fhn = H^1(X, \Omega^{[1,2\rangle }_X)^*$ and $\pair{,}$ is trivial, so that 
$\cz = H^1(X, \Omega^{[1,2\rangle }_X)^*$ and
 (\ref{comp-condition}) translates into the condition
 $
 \theta = \chi_0 \in  H^1(X, \Omega^{[1,2\rangle }_X)
 $.
 Thus, we can simplify the notation and write $\halfint_{\chi(z)}(\twcdo)$ instead of  $\halfintchitheta(\twcdo)$.

  Theorem \ref{equivalence-theorem} specializes to the following variant  of \cite{AChM}, Theorem 5.2.

\begin{thm}
Suppose that $\chi(z) \in H^1(X, \Omega^{[1,2\rangle }_X)((z))$ has regular singularity.
Then the functor
$$
\Sing: \halfintchitheta(\twcdo) \To \cD^{\chi_0}_X {\rm{-}mod}
$$
is an equivalence.
\end{thm}

When  $\cV$ is a CDO,  i.e. $\fh = 0$,
one obtains the  equivalence between half-integrable modules over a CDO  $\cV$
and the category of  $D$-modules over $X$.


\begin{thebibliography}{99}
\bibitem[AChM]{AChM} T.~Arakawa, D.~Chebotarov, F.~Malikov,
Algebras of twisted chiral differential operators and affine localization of $\fg$-modules,
arXiv:0810.4964, 
 to appear in Selecta Math.

\bibitem[BB1]{BB1} A.~Beilinson, J.~Bernstein, 
Localisation de $\frak{g}$-modules. (French)  C. R. Acad. Sci. Paris Se'r. I Math.
292  (1981),  no. 1, 15--18.


\bibitem[BB2]{BB2} A.~Beilinson, J.~Bernstein,
 A proof of Jantzen conjectures. I. M. Gelfand Seminar, 1--50,
 Adv. Soviet Math., 16, Part 1, Amer. Math. Soc., Providence, RI, 1993.



        \bibitem[Bo]{Borel} A.~Borel, ed.,
         Algebraic D-Modules, Perspectives in Mathematics, 2, Boston, MA: Academic Press, 1987.  ISBN 978-0-12-117740-9

	
	\bibitem[Bre]{Bre} P.~Bressler, The first Pontryagin class. Compos.
	Math.,   {\textbf{143}} (2007), 1127-1163
	
	\bibitem[Ch]{Ch} D.~Chebotarov, Classification of transitive vertex algebroids,
	arXiv:math.QA/1010.3385.

	\bibitem[F1]{F1} E.~Frenkel, Wakimoto modules, opers and the center at the
	critical level. Adv. Math. 195 (2005), no. 2, 297--404.

	
	\bibitem[F2]{F2} E.~Frenkel,
	 Lectures on the Langlands program and conformal field theory.
	 Frontiers in number theory, physics, and geometry. II, 387--533, Springer, Berlin,2007.
	 
	\bibitem[FBZ]{FBZ} E.~Frenkel, D.~Ben-Zvi, Vertex algebras and
	algebraic curves, 2nd edition, Mathematical Surveys and Monographs,
	v.58, AMS, 2004

	\bibitem[F3]{F3} E.~Frenkel, Langlands correspondence for loop groups,
	Cambridge University Press, 2007

	\bibitem[FF1]{FF1} B.~Feigin, E.~Frenkel, Representations of affine Kac-Moody
	algebras and bosonization, in: V.Knizhnik Memorial Volume,  L.Brink,
	D.Friedan, A.M.Polyakov (Eds.), 271-316, World Scientific,
	Singapore, 1990

	\bibitem[FF2]{FF2} B.~Feigin, E.~Frenkel, Affine Kac-Moody algebras at the
	critical level and Gelfand-Dikii algebras, in: Infinite Analysis,
	eds. A.Tsuchiya, T.Eguchi, M.Jimbo, {\it Adv. Series in Math. Phys.}
	{\bf 16} 197-215, Singapore, World Scientific, 1992


	

	\bibitem [GMS1]{GMS}  V.~Gorbounov, F.~Malikov, V.~Schechtman, Gerbes of chiral
	differential operators. II. Vertex algebroids, Invent. Math. 155
	(2004), no. 3, 605-680.

	\bibitem[GMS2]{GMSII} V.~Gorbounov, F.~Malikov, V.~Schechtman,
	 On chiral differential operators over homogeneous spaces.  Int. J. Math. Math. Sci.  26  (2001),  no.2, 83--106.

	\bibitem[HTT]{HTT}
	R.~Hotta, K.~Takeuchi, T.~Tanisaki, D-modules, perverse sheaves, and representation theory, Progress in Mathematics, 236, Boston, MA: Birkh\"auser Boston, 2008.  MR2357361, ISBN 978-0-8176-4363-8.
	
	\bibitem[K]{K} V.~Kac, Vertex algebras for beginners. University Lecture
	Series, 10. American Mathematical Society, Providence, RI, 1997.
	viii+141 pp. ISBN: 0-8218-0643-2

	\bibitem[KV1]{KV1} M.~Kapranov, E.~Vasserot,
	Vertex algebras and the formal loop space.  Publ. Math. Inst. Hautes \'Etudes Sci.   100  (2004), 209-269. 
	
	\bibitem[KV2]{KV2} M.~Kapranov, E.~Vasserot,
	 Formal loops IV: chiral differential operators,  preprint math.AG/0612371.
	 
 

	\bibitem[MSV]{MSV} F.~Malikov, V.~Schechtman, A.~Vaintrob,  Comm. in Math. Phys. {\bf 204} (1999), 439-473



\bibitem[W]{W} M.~Wakimoto, Fock representations of the affine Lie algebra
 $A\sp {(1)}\sb 1$. Comm. Math. Phys. {\bf 104} (1986), no. 4, 605--609.

	\bibitem[Z]{Zhu} Y.~Zhu, Modular invariance of characters of vertex operator
	algebras. J. Amer. Math. Soc. 9 (1996), no. 1, 237--302.
\end{thebibliography}
\end{document}